\documentclass[11pt,journal,onecolumn]{IEEEtran}   

\usepackage{amsmath,amssymb,bm,upgreek,cite,nicefrac,url}
\usepackage{dsfont}     %\mathds{1}
\usepackage[mathscr]{eucal}
\allowdisplaybreaks
\newcommand{\df}{:=}
\DeclareMathOperator{\Exp}{\mathbb{E}}
\DeclareMathOperator{\Prob}{\mathbb{P}}
\newcommand{\D}{\mathrm{d}}   % differential
\newcommand{\E}{\mathrm{e}}   % exponent
\newcommand{\RR}{\mathbb{R}}  % Real numbers
\newcommand{\Rd}{\mathbb{R}^d}
\newcommand{\NN}{\mathbb{N}}   % Natural numbers

\newcommand{\Ind}{\mathds{1}}        % indicator function
\newcommand{\Act}{\mathbb{U}}        % Action Set
\newcommand{\Uadm}{\mathfrak{U}}     % Admissible Controls
\newcommand{\Usm}{\mathfrak{U}_{\mathrm{SM}}}  % Stationary Markov controls
\newcommand{\Ussm}{\mathfrak{U}_{\mathrm{SSM}}}  % Stable stationary Markov controls

\newcommand{\transp}{^{\mathsf{T}}}  % transpose
\newcommand{\order}{{\mathscr{O}}}   % Order of
\newcommand{\grad}{\nabla}
\newcommand{\rc}{c}

\newcommand{\Bor}{\mathfrak{B}}     % Borel sets
\newcommand{\Cc}{\mathcal{C}}       % Set of continuous functions
\newcommand{\sF}{\mathfrak{F}}
\newcommand{\eom}{{\mathscr{G}}}   % Ergodic occupation measures
\newcommand{\cI}{\mathcal{I}}
\newcommand{\cK}{\mathcal{K}}
\newcommand{\Lg}{\mathcal{L}}       % Controlled extended generator
\newcommand{\Pm}{\mathfrak{P}}  % Probability measures
\newcommand{\cT}{\mathcal{T}}   % Modulated hitting times
\newcommand{\Lyap}{\mathcal{V}} % Lyapunov
\newcommand{\cX}{\mathcal{X}}   % Arbitrary space

\newcommand{\abs}[1]{\lvert#1\rvert}
\newcommand{\babs}[1]{\bigl\lvert#1\bigr\rvert}

\newcommand{\norm}[1]{\lVert#1\rVert}
\newcommand{\bnorm}[1]{\bigl\lVert#1\bigr\rVert}

\newtheorem{definition}{Definition}[section]
\newtheorem{lemma}{Lemma}[section]
\newtheorem{theorem}{Theorem}[section]

\newtheorem{remark}{Remark}[section]
\newtheorem{example}{Example}[section]
\newtheorem{assumption}{Assumption}[section]
\newtheorem{hypothesis}{Hypothesis}[section]
\newtheorem{corollary}{Corollary}[section]
\newcommand{\ttl}{\huge
Some new results on sample path optimality\\ in ergodic control of diffusions}
\begin{document}
\title{\ttl}

\author{\large Ari~Arapostathis\,$^*$\\
University of Texas at Austin, Austin, TX, USA
\thanks{$^\dagger$This work was supported in part by the Office of Naval Research
under the Electric Ship Research and Development Consortium, and in part by a grant from
the Army Research Office.}}

\markboth{Some new results on Sample Path Optimality in Ergodic Control}{}

\maketitle
\thispagestyle{empty}
%%%%%%%%%%%%%%%%%%%%%%%%%%%%%%%%%%%%%%%%%%%%%%%%%%%%%%%%%%%%%%%%%%%%%%%%%%%%%%%%
\begin{abstract}

We present some new results on sample path optimality
for the ergodic control problem of a class of non-degenerate diffusions
controlled through the drift.
The hypothesis most often used in the literature to ensure
the existence of an a.s.\ sample path optimal stationary Markov control
requires finite second moments of the first hitting
times $\uptau$ of bounded domains over all admissible controls.
We show that this can be considerably weakened: $\Exp[\uptau^2]$
may be replaced with $\Exp[\uptau\ln^+(\uptau)]$, thus reducing the
required rate of convergence of averages from polynomial to logarithmic.
A Foster--Lyapunov condition which guarantees this is also exhibited.
Moreover, we study a large class of models that are
neither uniformly stable, nor have a near-monotone running cost,
and we exhibit sufficient conditions for
the existence of a sample path optimal stationary Markov control.
\end{abstract}

\begin{IEEEkeywords}
ergodic control, sample path optimality, controlled diffusion,
subgeometric ergodicity,
empirical measures.
\end{IEEEkeywords}

\IEEEpeerreviewmaketitle
%%%%%%%%%%%%%%%%%%%%%%%%%%%%%%%%%%%%%%%%%%%%%%%%%%%%%%%%%%%%%%%%%%%%%%%%%%%%%%%%

\section{Introduction}

Sample path optimality in the ergodic control of diffusions has been studied
in \cite{Pra-00,Borkar,BorkGh-88a,book}.
For the analogous problem in Markov decision processes (MDP) we refer the reader to
\cite{Borkar-89,Cavazos-12,Cavazos-14,Cavazos-15,Lerma-99,Zhu-07,Hunt-05,Lasserre-99}.
In this paper we focus on non-degenerate diffusions with a compact action
space and controlled through the drift.

Consider a
controlled diffusion process $X = \{X_{t},\;t\ge0\}$
taking values in the $d$-dimensional Euclidean space $\RR^{d}$, and
governed by the It\^o stochastic differential equation
\begin{equation}\label{E-sde}
\D{X}_{t} \,=\, b(X_{t},U_{t})\,\D{t} + \upsigma(X_{t})\,\D{W}_{t}\,.
\end{equation}
All random processes in \eqref{E-sde} live in a complete
probability space $(\Omega, \sF,\Prob)$.
The process $W$ is a $d$-dimensional standard Wiener process independent
of the initial condition $X_{0}$.
The control process $U$ takes values in a compact, metrizable set $\Act$, and
$U_{t}(\omega)$ is jointly measurable in
$(t, \omega)\in[0,\infty)\times\Omega$.
Moreover, it is \emph{non-anticipative}:
for $s < t$, $W_{t} - W_{s}$ is independent of $\sF_{s}$ which is defined
as the the completion of $\sigma\{X_{0}, U_{r}, W_{r},\; r\le s\}$
relative to $(\sF,\Prob)$.
Such a process $U$ is called an \emph{admissible control},
and we let $\Uadm$ denote the set of all admissible controls.
We impose the usual assumptions on the data of the model to guarantee the existence and
uniqueness of solutions to \eqref{E-sde}.
These are described in Section~\ref{S-model}.

We adopt the relaxed control framework \cite[Section~2.3]{book}.
Thus, we extend the definition of
$U$ so that it takes values in $\Pm(\Act)$, the space of probability measures
on the Borel sets of $\Act$, and for a function $f\colon\Rd\times\Act$,
we let $f(\cdot\,,U_s) \equiv \int_{\Act} f(\cdot\,,u)\, U_s(\D{u})$.
For further details on relaxed controls see \cite[Section~2.3]{book}.
Since the set of optimal stationary controls, if nonempty, always contains
precise controls, the `relaxation' of the problem that allows
$\Pm(\Act)$-valued controls results in the same optimal value as the
original problem.
Thus, adopting relaxed controls
serves only to facilitate the analysis.

Let $\rc\colon \RR^{d} \times\Act\to\RR$ be a continuous function
bounded from below, which without loss of generality we assume is nonnegative,
referred to as the \emph{running cost}.
As is well known, the ergodic control problem, in its \emph{almost sure}
(or \emph{pathwise}) formulation,
seeks to a.s.\ minimize over all admissible controls $U$ the functional
\begin{equation}\label{E-ergcrit}
\limsup_{t\to\infty}\; \frac{1}{t}\int_{0}^{t} \rc(X_{s},U_{s})\,\D{s}\,.
\end{equation}
In equation \eqref{E-ergcrit}, $X_s$ denotes the process under the control $U$.

A weaker, \emph{average} formulation seeks to minimize
\begin{equation}\label{E-avgcrit}
\varrho_{U}\,\df\,\limsup_{t\to\infty}\;\frac{1}{t}\int_{0}^{t}
\Exp^{U}\bigl[\rc(X_{s},U_{s})\bigr]\,\D{s}\,.
\end{equation}
In this equation, $\Exp^{U}$ denotes the expectation operator associated
with the probability measure on the canonical
space of the process under the control $U$.
We define the \emph{optimal ergodic value} 
$\varrho_*\df\inf_{U\in\Uadm}\varrho_{U}$,
i.e., the infimum of \eqref{E-avgcrit} over all admissible controls.
We say that an admissible control $\Hat{U}$ is \emph{average cost optimal} if
$\varrho_{\Hat{U}}=\varrho_*$.
Also, an average cost optimal control $\Hat{U}$ is
called \emph{sample path
optimal}, or \emph{pathwise optimal}, if
\begin{equation}\label{E-pathopt}
\limsup_{t\to\infty}\; \frac{1}{t}\int_{0}^{t} \rc(X_{s},\Hat{U}_{s})\,\D{s}
\,\le\,
\liminf_{t\to\infty}\; \frac{1}{t}\int_{0}^{t} \rc(X_{s},U_{s})\,\D{s}
\end{equation}
a.s., for all $U\in\Uadm$.
It is evident that \eqref{E-pathopt} asserts a much stronger
optimality for $\Hat{U}$, viz., the most ``pessimistic'' pathwise cost 
under $\Hat{U}$ is no worse than the most ``optimistic'' pathwise cost under any
other admissible control.

For a control $U\in\Uadm$ we define the
process $\zeta^{U}_{t}$ of empirical measures as a
$\Pm(\RR^{d}\times\Act)$-valued process satisfying
\begin{equation}\label{E-emp}
\int_{\RR^{d}\times\Act} f\,\D{\zeta}^{U}_{t} \,=\, \frac{1}{t}\int_{0}^{t}
\int_{\Act} f(X_{s},u) U_{s}(\D{u})\,\D{s}\,,\quad t>0\,,
\end{equation}
for all $f\in\Cc_{b}(\RR^{d}\times\Act)$,
where $X$ denotes the solution of the diffusion in \eqref{E-sde} under
the control $U$.
Suppose for simplicity that the running cost is bounded.
It is then well known that a sufficient condition for the existence
of a \emph{pathwise
optimal} stationary Markov control is that
the family $\bigl\{\zeta^{U}_{t}\colon t\ge0\bigr\}$
is a.s.\ tight in $\Pm(\RR^{d}\times\Act)$ \cite[Theorem~3.4.7]{book}.
The hypothesis most often used in the literature
to guarantee the tightness of the family $\{\zeta^{U}_{t}\}$,
requires that the second moments of the first hitting time
of some bounded domain be
bounded uniformly over $X_0$ in a compact set, and all admissible controls $U\in\Uadm$.

We define
the family of operators $\Lg^{u}\colon\Cc^{2}(\RR^{d})\mapsto\Cc(\RR^{d})$ by
\begin{equation}\label{E-Lg}
\Lg^{u} f(x) \,=\, \frac{1}{2}\sum_{i,j} a^{ij}(x)
\frac{\partial^{2}f}{\partial x_{i} \partial x_{j}} (x)
+\sum_{i} b^{i} (x,u) \frac{\partial f}{\partial x_{i}} (x)
\end{equation}
for $u\in\Act$.
We refer to $\Lg^{u}$ as the \emph{controlled extended generator} of
the diffusion.
A Foster--Lyapunov condition which is sufficient for the uniform boundedness
of the second moments of hitting times is the following:
There exist a bounded domain $D\subset\Rd$,
and nonnegative functions $\Lyap_{1}$ and $\Lyap_{2}$ in
$\Cc^{2}(\Bar{D}^c)\cap\Cc(D^c)$
which satisfy
\begin{equation}\label{E-Lyap2nd}
\Lg^{u}\Lyap_{1} \,\le\, -1\,,\quad\text{and}\quad
\Lg^{u}\Lyap_{2} \,\le\, -\Lyap_{1}\qquad\text{on~} \Bar{D}^c\,,
\end{equation}
for all $u\in\Act$.
Here, $\Bar{D}$ and $D^c$ denote the closure and the complement of $D$,
respectively.
Let $\uptau(D)$ denote the first hitting time of the set $D$
for the process in \eqref{E-sde}, and
define $\Exp_x^U[\,\cdot\,]\df\Exp^U[\,\cdot\mid X_0=x]$,
with $x\in\Rd$.
The Foster--Lyapunov condition in \eqref{E-Lyap2nd}
implies that $\Exp_x^U[\uptau(D)\mid X_0=x]\le\Lyap_1(x)$,
and $\Exp_x^U\bigl[ \uptau(D)^2\bigr]\le2\Lyap_2(x)$
for all $U\in\Uadm$ and $x\in D^c$ \cite[Lemma~2.5.1]{book}.

The hypothesis that second moments of hitting times are bounded
can be considerably weakened. 
As we show in this paper, a sufficient condition for the
the family $\bigl\{\zeta^{U}_{t}\bigr\}$ to
be a.s.\ tight is
\begin{equation*}
\sup_{x\in K}\;\sup_{U\in\Uadm}\;
\Exp_{x}^{U} \bigl[\uptau(D)\,\ln^+(\uptau(D))\bigr]\,<\,\infty
\end{equation*}
for any compact set $K\subset D^c$,
where $\ln^+$ denotes the positive part of the natural logarithm.
We also show that the second inequality in \eqref{E-Lyap2nd}
may be replaced  by
 $\Lg^{u}\Lyap_{2} \le -\ln^+(\Lyap_{1})$.

\subsection{The model}\label{S-model}
The following assumptions are in effect throughout the paper.

\begin{description}
\item[(A1)]
\textit{Local Lipschitz continuity:\/}
The functions
\begin{equation*}
b=\bigl[b^{1},\dotsc,b^{d}\bigr]\transp \colon\RR^{d} \times\Act\mapsto\RR^{d}
\end{equation*}
and $\upsigma=\bigl[\upsigma^{ij}\bigr]\colon\RR^{d}\mapsto\RR^{d\times d}$
are locally Lipschitz in $x$ with a Lipschitz constant $\kappa_{R}$ depending on
$R>0$.
In other words, if $B_{R}$ denotes the open ball of radius $R$ centered
at the origin in $\RR^{d}$, then
for all $x,y\in B_{R}$ and $u\in\Act$,
\begin{equation*}
\abs{b(x,u) - b(y,u)} + \norm{\upsigma(x) - \upsigma(y)}
\,\le \,\kappa_{R}\abs{x-y}\,,
\end{equation*}
where $\norm{\upsigma}^{2}\df
\mathrm{trace}\left(\upsigma\upsigma\transp\right)$.
We also assume that $b$ is continuous in $(x,u)$.

\item[(A2)]
\textit{Affine growth condition:\/}
$b$ and $\upsigma$ satisfy a global growth condition of the form
\begin{equation*}
\abs{b(x,u)}^{2}+ \norm{\upsigma(x)}^{2}\,\le\, \kappa_{1}
\bigl(1 + \abs{x}^{2}\bigr)
\end{equation*}
for all $(x,u)\in\RR^{d}\times\Act$.

\item[(A3)]
\textit{Local non-degeneracy:\/}
Let $a \df \upsigma\,\upsigma\transp$.
For each $R>0$, we have
\begin{equation*}
\sum_{i,j=1}^{d} a^{ij}(x)\xi_{i}\xi_{j}
\,\ge\,\kappa^{-1}_{R}\abs{\xi}^{2} \qquad\forall x\in B_{R}\,,
\end{equation*}
for all $\xi=(\xi_{1},\dotsc,\xi_{d})\in\RR^{d}$.

\item[(A4)]
The running
cost function $\rc\colon\RR^{d}\times\Act\to\RR_{+}$
is continuous.
\end{description}

The conditions in (A1)--(A3) are standard assumptions on the drift $b$
and the diffusion matrix $\upsigma$
to guarantee existence and uniqueness of solutions to \eqref{E-sde}.
The running cost $\rc$ is usually assumed to be locally Lipschitz in
its first argument.
However, this is only used to obtain regular solutions to the Hamilton--Jacobi--Bellman
equation, something which does not concern us in this paper.

\subsection{Markov controls}
Of fundamental importance in the study of functionals of $X$ is
It\^o's formula.
For $f\in\Cc^{2}(\RR^{d})$ and with $\Lg^{u}$ as defined in \eqref{E-Lg},
it holds that
\begin{equation}\label{E-Ito}
f(X_{t}) \,=\, f(X_{0}) + \int_{0}^{t}\Lg^{U_{s}} f(X_{s})\,\D{s}
+ M_{t}\,,\quad\text{a.s.},
\end{equation}
where
\begin{equation*}
M_{t} \,\df\, \int_{0}^{t}\bigl\langle\nabla f(X_{s}),
\upsigma(X_{s})\,\D{W}_{s}\bigr\rangle
\end{equation*}
is a local martingale.

Recall that a control is called \emph{Markov}  if
$U_{t} = v(t,X_{t})$ for a measurable map
$v \colon\RR\times\RR^{d}\mapsto\Pm(\Act)$,
and it is called \emph{stationary Markov} if $v$ does not depend on
$t$, i.e., $v \colon\RR^{d}\mapsto\Pm(\Act)$.
It is well known that under Assumptions (A1)--(A3),
under any Markov control $v$, the stochastic differential equation in
\eqref{E-sde} has a unique strong solution \cite{Gyongy-96}.

Let $\Usm$ denote the set of stationary Markov controls.
Under any $v\in\Usm$, the process $X$ is strong Markov,
and we denote its transition function by $P^{v}_{t}(x,\cdot)$.
It also follows from the work of \cite{Bogachev-01,Stannat-99} that under any
$v\in\Usm$, the transition probabilities of $X$
have densities which are locally H\"older continuous.
Thus $\Lg^{v}$ defined by
\begin{equation*}
\Lg^{v} f(x) \,=\, \sum_{i,j} a^{ij}(x)
\frac{\partial^{2}f}{\partial x_{i} \partial x_{j}} (x)
+\sum_{i} b^{i} (x,v(x)) \frac{\partial f}{\partial x_{i}} (x)
\end{equation*}
for $v\in\Usm$ and $f\in\Cc^{2}(\RR^{d})$,
is the generator of a strongly-continuous
semigroup on $\Cc_{b}(\RR^{d})$, which is strong Feller,
i.e., if $f$ is bounded measurable function $f$, and $t>0$, then
the map
$x\mapsto \Exp[ f(X_t) \mid X_0=x]$ is continuous.
We let $\Prob_{x}^{v}$ denote the probability measure and
$\Exp_{x}^{v}$ the expectation operator on the canonical space of the
process under the control $v\in\Usm$, conditioned on the
process $X$ starting from $x\in\RR^{d}$ at $t=0$.
The operator $\Lg^{U}$ for $U\in\Uadm$ is analogously defined.

In the next section, 
we summarize the notation used in the paper.

\subsection{Notation}

The standard Euclidean norm in $\RR^{d}$ is denoted by $\abs{\,\cdot\,}$,
and $\RR_{+}$ stands for the set of non-negative real numbers.
The closure, the complement, and the indicator function
of a set $A\subset\Rd$ are denoted
by $\Bar{A}$, $A^{c}$, and $\Ind_A$, respectively.
We denote by $\uptau(A)$ the first hitting time of $A\subset\RR^{d}$, i.e.,
\begin{equation*}
\uptau(A) \;\df\; \inf\;\{t\ge0\;\colon\, X_{t}\in A\}\,.
\end{equation*}
The open ball in $\RR^{d}$ centered at the origin, of radius $R>0$,
is denoted by $B_{R}$.

The Borel $\sigma$-field of a topological space $E$ is denoted by $\Bor(E)$,
and $\Pm(E)$ stands for the set of
probability measures on $\Bor(E)$.
The space $\Pm(E)$ is always viewed as endowed with the topology of
weak convergence of probability measures (the Prohorov topology).

We introduce the following notation for spaces of real-valued functions on
a domain $D\subset\RR^{d}$.
The space $\Cc^{k}(D)$
refers to the class of all functions whose partial
derivatives up to order $k$ exist and are continuous,
and $\Cc_{b}^{k}(\RR^{d})$
is the subspace of $\Cc^{k}(\RR^{d})$
consisting of those functions whose derivatives up to order $k$ are bounded.

We also need the following definition.

\begin{definition}
A function $h\colon\cX\to\RR$, where $\cX$ is a locally compact
space, is called \emph{inf-compact}
on a set $A\subset \cX$ if the set
$\Bar{A}\cap\bigl\{y\,\colon\,h(y)\le \beta\bigr\}$
is compact (or empty) in $\cX$ for all $\beta\in\RR$.
When this property
holds for $A=\cX$, then we simply say that $h$ is inf-compact.
\end{definition}

\section{Sample Path Optimality for Uniformly Stable Diffusions}
\label{S2}

Recall that a control $v\in\Usm$ is called \emph{stable}
if the associated diffusion is positive recurrent.
We let $\Ussm$ denote the set of stable controls,
and let $\mu_{v}$ stand for the unique invariant probability
measure of the diffusion process under the control $v\in\Ussm$.
For $v\in\Ussm$, we define the \emph{ergodic occupation measure}
$\uppi_{v}\in\Pm(\Rd\times\Act)$ by
$$\uppi_{v}(\D{x},\D{u})\,\df\,\mu_{v}(\D{x})\,v(\D{u}\mid x)\,,$$
and let $\eom$ denote the set of ergodic occupation measures
(see Section~3.2.1 in \cite{book}).
It is well known that that
$\eom$ is a closed and convex subset of $\Pm(\RR^{d}\times\Act)$
 \cite[Lemma 3.2.3]{book}.
The controlled diffusion in \eqref{E-sde} is called \emph{uniformly stable}
if the set of ergodic occupation measures is compact.
It is well known that the following condition is sufficient for
uniform stability:
$\uptau(D)$ is uniformly integrable with respect to
$\{\Prob^v_x \colon v\in\Usm\}$
for some bounded open domain $D$, and some $x\in \Bar{D}^c$.
Conversely, if the controlled diffusion is uniformly stable
then the uniform integrability property holds with respect to
the larger family $\{\Prob^U_x \colon U\in\Uadm\}$,
for any bounded open domain $D$ and any $x\in \Bar{D}^c$ over 
\cite[Lemma~3.3.4]{book}.

It is well known (see \cite[p.~19]{Meyer-66}) that uniform integrability of $\uptau(D)$
is equivalent to the existence
of some nonnegative $\phi\in\Cc(\RR_+)$ satisfying $\phi(z)\to\infty$
as $z\to\infty$, and
\begin{equation}\label{E-unif}
\sup_{U\in\Uadm}\;\Exp_x^U\bigl[\uptau(D)\phi\bigl(\uptau(D)\bigr)\bigr]\,<\,\infty\,.
\end{equation}
A question posed in \cite[Remark~5.10\,(i)]{Survey} in the context of
countable MDPs was whether uniform integrability of the
hitting times is sufficient for pathwise optimality.
A counterexample to this appeared very recently in \cite{Cavazos-14}.
Note though that the optimal control problem in this example
is equivalent to minimizing the ``$\liminf$'' for a function which is
bounded above.
Even though this counterexample cannot be adapted for diffusions due to the
nature of the transition probabilities, it suggests 
that \eqref{E-unif} may not be sufficient for pathwise optimality.
It turns out, however, that if \eqref{E-unif} holds with $\phi=\ln^+$,  uniformly
over all $x$ in compact subsets of $D^c$, then the family of empirical
measures
$\{\zeta^{U}_{t}\colon\;t\ge0\}$ defined in \eqref{E-emp}
is a.s.\ tight for any $U\in\Uadm$.

As seen from the proof of \cite[Theorem~3.4.11]{book} the hypothesis
of uniform boundedness of second moments of hitting times is used in
applying Kolmogorov's strong law of large numbers for martingale difference
sequences (for more details see \cite[Lemma~2.12]{Borkar-02}).
But this can be accomplished under weaker hypotheses.
Let $\{Y_n\}$ be a sequence of
integrable random variables
defined on a complete probability space $(\Omega,\sF,\Prob)$,
and define $\sF_{n}\df \sigma(Y_{1},Y_{2},\dotsc,Y_{n})$, and
\begin{equation*}
S_{n} \,\df\, \sum_{k=1}^{n}\bigl(Y_{k}-
\Exp\left[Y_{k}\mid \sF_{k-1}\right]\bigr)\,, \quad n\in\NN\,.
\end{equation*}
Fairly recently, Stoica has established in \cite[Theorem~2]{Stoica-11} that if
$\sup_{n\in\NN}\, \Exp\bigl[\abs{Y_n}\,\ln^+\abs{Y_n}\bigr]<\infty$, then
it holds that
\begin{equation*}
\frac{S_{n}}{n}\;\xrightarrow[n\to\infty]{\text{a.s.}}\;0\,.
\end{equation*}
This result is sharp.
As shown in \cite{Elton-81}, given any centered random variable $Z$
satisfying $\Exp\bigl[\abs{Z}\bigr]<\infty$,
and $\Exp\bigl[\abs{Z}\,\ln^+\abs{Z}\bigr]=\infty$, there exists
a martingale difference sequence $\{Z_n\}$ of identically
distributed random variables,  having
the same law as $Z$, such that
$\frac{Z_1+\dotsb+Z_n}{n}$ diverges almost surely as $n\to\infty$.

We present the following theorem.

\begin{theorem}\label{T2.1}
Suppose that for some bounded domain $D\subset\Rd$,
and any compact set $K\subset D^c$, it holds that
\begin{equation*}
\sup_{x\in K}\;\sup_{U\in\Uadm}\;
\Exp_x^U\bigl[\uptau(D)\ln^+\bigl(\uptau(D)\bigr)\bigr]\,<\,\infty\,.
\end{equation*}
Then $\{\zeta^{U}_{t}\colon\;t\ge0\}$ is a.s.\ tight for any $U\in\Uadm$.
\end{theorem}

\begin{IEEEproof}
We follow the proof of \cite[Theorem~3.4.11]{book}.
Without loss of generality we may assume $D$ is an open ball.
Let $\Tilde{D}$ be a ball containing $\Bar{D}$.
Let $\Hat{\uptau}_{0}=0$, and for $k=0,1,\dotsc$ define inductively
an increasing sequence of stopping times by
\begin{equation*}\label{E2.6.11}
\begin{split}
\Hat{\uptau}_{2k+1}
&\df\inf\;\{t>\Hat{\uptau}_{2k}: X_{t}\in \Tilde{D}^{c}\}\,,\\[3pt]
\Hat{\uptau}_{2k+2}
&\df \inf\; \{t>\Hat{\uptau}_{2k+1}: X_{t}\in D\}\,.
\end{split}
\end{equation*}
The assumption of the theorem implies of course that
$\uptau(D)$ satisfies \eqref{E-unif}, or, equivalently that
it is uniformly integrable with respect to
$\{\Prob^v_x \colon v\in\Usm\}$.
Thus, by Lemma~3.3.4 in \cite{book} there exist inf-compact functions
$\Phi\in\Cc^{2}(\Rd)$ and
$h\colon\Rd\times\Act\to\RR_+$, and a constant $k_0$, satisfying
\begin{equation*}
\Lg^u\Phi(x) \;\le\; k_0 - h(x,u) \qquad\forall (x,u)\,\in\Rd\times\Act\,.
\end{equation*}
This implies that
\begin{equation*}
\sup_{x\in\partial{\Tilde{D}}}\;\sup_{U\in\Uadm}\;
\Exp_{x}^{U}\left[\int_{0}^{\uptau(D)}
h(X_{t},U_t)\,\D{t}\right]<\infty\,.
\end{equation*}
Since $h$ is inf-compact, and also
$\displaystyle\inf_{x\in\partial{\Tilde{D}}}\;\inf_{U\in\Uadm}\;
\Exp_{x}^{U}\left[\uptau(D)\right]>0$, it follows
that set  of probability measures $\nu$ defined by
\begin{equation}\label{ET2.1A}
\int_{\RR^{d}} f\,\D{\nu} =
\frac{\Exp_{X_{0}}^{U}\left[\int_{0}^{\Hat{\uptau}_{2}} f(X_{t})\,\D{t}\right]}
{\Exp_{X_{0}}^{U}\bigl[\Hat{\uptau}_{2}\bigr]}
\qquad\forall f\in\Cc_{b}(\RR^{d})\,,
\end{equation}
for $U\in\Uadm$, and with the law of $X_{0}$
supported on $\partial{D}$, is tight.
By \cite[Theorem~2.6.1\,(b)]{book}, we have
$$\displaystyle\sup_{x\in\partial{D}}\;\inf_{U\in\Uadm}\;
\Exp_{x}^{U}\left[\uptau^k(\Tilde{D}^{c})\right]<\infty$$ for any $k\in\NN$.
This, together with the hypothesis of the theorem, imply that
\begin{equation*}
\sup_{x\in D}\;\sup_{U\in\Uadm}\;
\Exp_x^U\bigl[\Hat{\uptau}_{2}
\ln^+\bigl(\Hat{\uptau}_{2}\bigr)\bigr]\,<\,\infty\,.
\end{equation*}
Therefore, we can use the strong law of large numbers for a martingale
difference sequence to deduce that,
for all $f\in\Cc_{b}(\RR^{d})$,
\begin{equation}\label{ET2.1B}
\frac{1}{n}\int_{0}^{\Hat{\uptau}_{2n}} f(X_{t})\,\D{t}\\[5pt]
-\frac{1}{n}\sum_{m=0}^{n-1}
\Exp_{X_{0}}^{U}\left[\int_{\Hat{\uptau}_{2m}}^{\Hat{\uptau}_{2m+2}} 
f(X_{t})\,\D{t}\Bigm|\sF_{\Hat{\uptau}_{2m}}\right]
\xrightarrow[n\to\infty]{\text{a.s.}} 0
\end{equation}
for all $f\in\Cc_{b}(\RR^{d})$.
We then proceed as in the proof of \cite[Theorem~3.4.11]{book} to show,
using \eqref{ET2.1B} and the tightness of the measures $\nu$ defined in \eqref{ET2.1A},
that for any $\epsilon>0$ there exists
$N(\epsilon)\in\NN$ such that if $f_\ell\colon\Rd\to[0,1]$ is any
continuous function which vanishes on $B_\ell$ and is equal to $1$ on $B^c_{\ell+1}$,
then
\begin{equation*}
\limsup_{n\to\infty}\;
\frac{1}{\Hat{\uptau}_{2n}}\int_{0}^{\Hat{\uptau}_{2n}} f_{\ell}(X_{t})\,\D{t}
\;\le\;\epsilon\quad\text{a.s.,}\quad\forall\,\ell\ge N(\epsilon)\,.
\end{equation*}
Then a standard approximation argument using
the number of cycles completed at time $t$ defined as
$\varkappa(t) = \max\; \{k: t> \Hat{\uptau}_{2k}\}$ shows that
\begin{equation*}
\limsup_{t}\;
\frac{1}{t}\int_{0}^{t} f_{\ell}(X_{t})\,\D{t}
\;\le\;\epsilon\quad\text{a.s.,}\quad\forall\,\ell\ge N(\epsilon)\,,
\end{equation*}
which implies that the family of empirical measures
$\{\zeta^{U}_{t}\}$ is a.s.\ tight.
This completes the proof.
\end{IEEEproof}

We continue with the following lemma.

\begin{lemma}\label{L2.1}
Let $D$ be a bounded $\Cc^{2}$ domain.
If there exist nonnegative functions $\varphi_{1}$ and
$\varphi_{2}$ in $\Cc^{2}(\Bar{D}^{c})\cap\Cc(D^{c})$
satisfying, for some $v\in\Usm$,
\begin{equation}\label{EL2.1a}
\begin{split}
\Lg^{v}\varphi_{1}(x)&\,\le\, -1\,,\\[5pt]
\quad\Lg^{v}\varphi_{2}(x)&\,\le\, -\ln^+\bigl(\varphi_{1}(x)\bigr)\,,
\end{split}
\end{equation}
for all $x\in \Bar{D}^c$, then
\begin{equation}\label{EL2.1b}
\begin{split}
\Exp_{x}^{v}\bigl[\uptau(D)\bigr]&\,\le\, \varphi_{1}(x)\,,\\[5pt]
\Exp_{x}^{v}\bigl[\uptau(D)\ln^+\bigl(\uptau(D)\bigr) \bigr]&\,\le\,
2\varphi_{1}(x)+\varphi_{2}(x)\,,
\end{split}
\end{equation}
for all $x\in D^{c}$.
\end{lemma}

\begin{IEEEproof}
Without loss of generality, since $\varphi_1$ and $\varphi_2$ can always
be scaled, we may assume that $\varphi_1\ge1$ on $D^c$.
The first inequality in \eqref{EL2.1b} is standard.

It is simple to verify that for any $T>0$ we have
\begin{equation}\label{EL2.1c}
T\,\ln^+T \,\le\, \int_0^T \ln^+(T-t)\,\D{t} + T\,.
\end{equation}
Let $\Hat{\uptau}_{R} \df \uptau(D)\wedge \uptau(B^c_R)$, $R>0$.
Using Dynkin's formula, we obtain
\begin{equation}\label{EL2.1d}
\Exp_{x}^{v} \left[ \int_{0}^{\Hat{\uptau}_{R}}
\ln\bigl(\varphi_{1}(X_{t})\bigr)\,\D{t}\right]
\,\le\, \varphi_{2}(x)\,.
\end{equation}
Conditioning at  $t\wedge\Hat{\uptau}_{R}$, and using optional
sampling and Jensen's inequality, we obtain
\begin{align}\label{EL2.1e}
\Exp_{x}^{v}\biggl[\int_0^{\Hat{\uptau}_{R}} \ln^+(\Hat{\uptau}_{R}-t)\,\D{t}\biggr]
&\,\le\, \Exp_{x}^{v}\biggl[\int_{0}^{\infty} \ln(\Hat{\uptau}_{R}-t+1)\,
\Ind_{\{t<\Hat{\uptau}_{R}\}}\,\D{t} \biggr]\nonumber \\
&\,=\, \Exp_{x}^{v}\biggl[ \int_{0}^{\infty} \Exp_{x}^{v}
\Bigl[\ln(\Hat{\uptau}_{R}-t+1)\,
\Ind_{\{t<\Hat{\uptau}_{R}\}} \Bigm|
\sF^{X}_{t\wedge\Hat{\uptau}_{R}}\Bigr] \,\D{t} \biggr]\nonumber \\
&\,\le\, \Exp_{x}^{v}\biggl[ \int_{0}^{\infty} 
\Ind_{\{t\wedge\Hat{\uptau}_{R}<\Hat{\uptau}_{R}\}}
\ln\bigl(\Exp^{v}_{X_{t\wedge\Hat{\uptau}_{R}}}
\bigl[\Hat{\uptau}_{R}+1\bigr]\bigr) \,\D{t} \biggr]\nonumber \\
&\,\le\, \Exp_{x}^{v}\biggl[\int_{0}^{\infty}
\ln\bigl(\varphi_{1}(X_{t\wedge\Hat{\uptau}_{R}})+1\bigr)\,
\Ind_{\{t<\Hat{\uptau}_{R}\}}\,\D{t} \biggr]\nonumber \\
&\,\le\, \Exp_{x}^{v}\left[\int_{0}^{\Hat{\uptau}_{R}}
\bigl(\ln\bigl(\varphi_{1}(X_{t})\bigr)+1\bigr)\,\D{t}\right]\,,
\end{align}
where we use the inequality $\ln(z)+1\ge \ln(z + 1)$ for all
$z\in[1,\infty)$.
Combining \eqref{EL2.1c}--\eqref{EL2.1e}, and letting $R\to\infty$, we obtain
the second inequality in \eqref{EL2.1b}.
\end{IEEEproof}

\begin{remark}\label{R2.1}
It is evident from the proof of Lemma~\ref{L2.1} that if
we replace $v$ in \eqref{EL2.1a} by some $U\in\Uadm$, then
\eqref{EL2.1b} holds for the process controlled under $U$.
Likewise, if \eqref{EL2.1a} holds for all $u\in\Act$, then
\eqref{EL2.1b} holds uniformly over $U\in\Uadm$.
\end{remark}

Consider the following hypothesis.

\begin{hypothesis}\label{H01}
There exist a bounded domain $D\subset\Rd$, and positive
functions $\Lyap_{1}$ and $\Lyap_{2}$,
with $\Lyap_{i}\in\Cc^{2}(\RR^{d})$, $i=1,2$, such that
\begin{subequations}
\begin{align}
\Lg^{u}\Lyap_{1}(x) &\,\le\, - 1\,,\label{EH01a}\\[5pt]
\Lg^{u}\Lyap_{2}(x) &\,\le\, -\ln^+\bigl(\Lyap_{1}(x)\bigr)\label{EH01b}
\end{align}
\end{subequations}
for all $x\in D^c$ and $u\in\Act$.
\end{hypothesis}

We have the following corollary.

\begin{corollary}\label{C2.1}
Let Hypothesis~\ref{H01} hold, and suppose that
there exists $\Hat{v}\in\Ussm$ such that $\varrho_{\Hat{v}}$ is finite.
Then
\begin{itemize}
\item[(a)]
There exists a stationary Markov control which is average cost optimal.
\item[(b)]
Every average cost optimal stationary Markov control is pathwise optimal.
\end{itemize}
\end{corollary}

\begin{IEEEproof}
As shown in \cite[p.~65]{book} for any bounded domain $D$ it holds
that $\Exp_x^U\bigl[\uptau(D)]\to\infty$ as $\abs{x}\to\infty$.
Since by It\^o's formula and \eqref{EH01a} we have
$\Lyap_{1}(x)\ge\Exp_x^U\bigl[\uptau(D)]$ for all $U\in\Uadm$, it
follows that $\Lyap_1$ is inf-compact.
Therefore, \eqref{EH01b} implies that the set of ergodic
occupation measures $\eom$ is compact (see Lemma~3.3.4\,(iv) in \cite{book}).
Since $\rc$ is bounded below, part~(a) follows by \cite[Theorem~3.4.5]{book}.

Part~(b) follows by Theorem~\ref{T2.1}, Lemma~\ref{L2.1},
and Remark~\ref{R2.1}.
\end{IEEEproof}

The pair of Lyapunov equations in
Hypothesis~\ref{H01} may be replaced by a single equation,
which in many cases it might be easier to verify (see Example~\ref{E2.1}
later in this section).
Consider the following hypothesis.

\begin{hypothesis}\label{H02}
There exist an inf-compact function $\Lyap\in\Cc^{2}(\RR^{d})$,
with $\Lyap\ge1$,
and a constant $C$ such that
\begin{equation}\label{E-Lyap3}
\Lg^u\Lyap(x) \,\le\,
C - \ln\bigl(\Lyap(x)\bigr)\qquad\forall\,(x,u)\in\Rd\times\Act\,.
\end{equation}
\end{hypothesis}

Let $D$ be an open ball such that $\ln\bigl(\Lyap(x)\bigr)\ge 2C+1$
for $x\in D^c$.
Then, if we define $\Lyap_{1}\df\Lyap$ and $\Lyap_{2}\df 2\Lyap$,
it is clear that \eqref{E-Lyap3} implies \eqref{EH01a}--\eqref{EH01b}.
Thus Hypothesis~\ref{H02} implies Hypothesis~\ref{H01}.
The Foster--Lyapunov condition in \eqref{E-Lyap3} results in
subgeometric ergodicity for the process
\cite[Theorems~3.2 and 3.4]{Douc-09}.
A similar estimate as in Lemma~\ref{L2.1} can be derived directly
from \eqref{E-Lyap3},
albeit we have to restrict to stationary controls.
Without loss of generality, let $\phi(y)\df 1+\ln(y)$, $y\in[1,\infty)$
and select $D$ and a constant $C$ such that
\begin{equation}\label{E-FLlog}
\Lg^u\Lyap(x) \,\le\,
C\,\Ind_D(x) - \phi\bigl(\Lyap(x)\bigr)\qquad\forall\,(x,u)\in\Rd\times\Act\,.
\end{equation}
Let $\cI\colon[1,\infty)\to\RR_+$ denote the `shifted' logarithmic integral
\begin{equation*}
\cI(z)\,\df\, \int_{1}^z \frac{\D{s}}{1 +\ln(s)}\,.
\end{equation*}
Combining the identity $\bigl(\cI^{-1})'(z) = \phi\bigl(\cI^{-1}(z)\bigr)$
and \cite[Theorem~4.1]{Douc-09} we obtain
\begin{equation*}
\Exp_{x}^{v}\bigl[\cI^{-1}\bigl(\uptau(D)\bigr) \bigr]\,\le\, \Lyap(x)-1
\qquad\forall\,x\in D^c\,,\quad\forall v\in\Usm\,.
\end{equation*}
Since $\cI^{-1}(z)$ grows as $\order(z\ln(z))$, the estimate of
Theorem~\ref{T2.1} follows.

The Foster--Lyapunov condition in
\eqref{E-FLlog} implies that,
under any stationary control, the process is ergodic
at a logarithmic rate.
Indeed, applying \cite[Theorem~3.2]{Douc-09} with $\varPsi_1$ the identity
function, and $\varPsi_2=1$,
it follows by (3.5) in \cite{Douc-09} that there exists a positive
constant $C_0$,
such that
\begin{equation*}
\bnorm{P^{v}_{t}(x,\cdot\,)- \mu_v(\,\cdot\,)}_{\text{TV}}
\;\le\; \frac{C_0}{\ln(t+1)}\,\Lyap(x)\quad\forall\, (t,x)\in(0,\infty)\times\Rd\,,
\end{equation*}
and for all $v\in\Usm$,
where $\bnorm{\,\cdot\,}_{\text{TV}}$ stands for the total variation norm.
See also \cite[p.~1364]{Douc-04} for the
corresponding results for Markov chains.
For some other recent developments in this topic see
\cite{Douc-04,Bakry-08, Douc-09,Fort-05,Hairer-16}.

The following is an example of a uniformly stable controlled diffusion which satisfies
$\Exp^v_x[\uptau(D)^p]=\infty$ for all $p>1$, and all $v\in\Usm$, but
$\Exp_{x}^{v}\bigl[\uptau(D)\ln^+\bigl(\uptau(D)\bigr) \bigr]<\infty$.

\begin{example}\label{E2.1}
Consider the one-dimensional controlled diffusion
\begin{equation}\label{E-sde2}
\D{X}_{t} \,=\,
- \frac{X_{t}\Bigl(1+\bigl[\ln(2+X_t^2)\bigr]^{-\nicefrac{1}{2}}\,U_{t}\Bigr)}
{2+ X_t^2}\,\D{t}
+ \sqrt{2}\,\D{W}_{t}\,,
\end{equation}
with $\Act=[1,2]$.
We apply Theorem~3\,(c) in \cite{Balaji-00-passage} with
$A(x)\df 2$, $B(x)\df 2$ and
\begin{equation*}
C(x) \;\df\; -\frac{x^2\Bigl(1+\bigl[\ln(2+x^2)\bigr]^{-\nicefrac{1}{2}}\,v(x)\Bigr)}
{2+ x^2}\,.
\end{equation*}
It then follows
that for any $v\in\Usm$, $r>0$, $\abs{x}>r$, and $p>1$, we have
$\Exp^v_x[\uptau(B_r)^p]=\infty$.
Here, $B_r=[-r,r]$.
Thus the hitting times do not have any moments larger than $1$.

Let $\Lyap(x) \df 1+x^2\bigl[\ln(2+x^2)\bigr]^2$.
A straightforward calculation shows that for some large
enough constant $\kappa>0$
we have
\begin{equation}\label{EE2.1A}
\Lg^u\Lyap(x) \,\le\, \kappa - \frac{3}{2}\,\bigl[\ln(2+x^2)\bigr]^{\nicefrac{3}{2}}
\qquad\forall\,(x,u)\in\Rd\times\Act\,.
\end{equation}
This of course can be written as
\begin{equation*}
\Lg^u\Lyap(x) \,\le\, C - \ln\bigl(\Lyap(x)\bigr)
\qquad\forall\,(x,u)\in\Rd\times\Act\,,
\end{equation*}
for some large enough constant $C$.
Therefore, Hypothesis~\ref{H02} is satisfied.
Let $h(x)\df\bigl[\ln(2+x^2)\bigr]^{\nicefrac{3}{2}}$.
By \eqref{EE2.1A} and It\^o's formula we obtain
$\int_{\RR} h(x)\,\mu_v(\D{x})\le\frac{2}{3}\kappa$ for all $v\in\Usm$.
In addition, \eqref{EE2.1A} implies that the set of ergodic occupation measures
is compact.
It then follows by Corollary~\ref{C2.1}\,(a),
that if $\rc$ is a running cost which is bounded below
and such that $x\mapsto\max_{u\in\Act}\,\abs{\rc(x,u)}$ does
not grow faster than $\bigl(\ln^+\abs{x}\bigr)^{\nicefrac{3}{2}}$,
then there exists $v^*\in\Ussm$ which is average cost optimal for
the diffusion in \eqref{E-sde2} \cite[Theorem~3.4.5]{book}.
Moreover, every average cost optimal stationary Markov control
is pathwise optimal by Corollary~\ref{C2.1}\,(b).
\end{example}

\begin{remark}
If the running cost $\rc$ is not bounded below, one may obtain the
same results by considering the modulated hitting times
\begin{equation*}
\cT(D) \,\df\, \int_{0}^{\uptau(D)} \max_{u\in\Act}\,\abs{\rc(X_t,u)}\,\D{t}\,.
\end{equation*}
Theorem~\ref{T2.1} then still holds if we replace $\uptau$ with $\cT$ in
the hypothesis.
Analogously, if we replace \eqref{EH01a} by
\begin{equation*}
\Lg^{u}\Lyap_{1}(x) \,\le\, - g(x)\qquad \forall\,(x,u)\in\Bar{D}^c\times\Act\,,
\end{equation*}
where $g\ge1$
is some function  that grows at least as fast as the map
$x\mapsto\max_{u\in\Act}\,\abs{\rc(x,u)}$,
then Corollary~\ref{C2.1} remains valid for a running cost $\rc$ which is not
necessarily bounded below.
\end{remark}

\section{Sample path optimality for a general class of diffusions}

The running cost is called \emph{near-monotone} if there exists some
$\epsilon>0$ such that the level set $\{(x,u)\in\Rd\times\Act\colon
\rc(x,u)\le \varrho_*+\epsilon\}$ is compact.
It is well known that if the running cost is near-monotone then
every average cost optimal stationary Markov control is necessarily pathwise optimal
\cite[Theorem~3.4.7]{book}.
In this Section we consider a general class of ergodic control problems
for which the diffusion is not uniformly stable, and the running cost is not
near monotone.
They are characterized by the property that the running cost is near monotone
when restricted to some subset $\cK\subset\Rd$, while a suitable Foster--Lyapunov
condition holds on $\cK^c$.
Models of this nature appear, for example, as the limiting diffusions of
multiclass queueing networks in the Halfin--Whitt regime \cite{ABP14}.

Consider the following assumption.
\begin{assumption}\label{A3.1}
For some open subset $\cK\subset\Rd$, the following hold.
\begin{itemize}
\item[(i)]
There exists $\Hat{U}\in\Uadm$ such that $\varrho_{\Hat{U}}<\infty$.
\item[(ii)]
The level sets $\{(x,u)\in\cK\times\Act\colon \rc(x,u)\le\gamma\}$ are
bounded (or are empty) for all $\gamma\in(0,\infty)$\,.

\item[(iii)]
For some smooth increasing concave function
$\phi\colon\RR_+\to\RR_+$, satisfying $\phi(z)\to\infty$ as $z\to\infty$,
some positive inf-compact function $\Lyap\in\Cc^{2}(\RR^{d})$,
and some positive constants $\kappa_i$, $i=0,1,2$, it holds that
\begin{equation}\label{E-KA1}
\begin{split}
\Lg^{u}\Lyap(x) &\,\le\,\kappa_0 - \kappa_1\,\phi\bigl(\Lyap(x)\bigr)\qquad
\forall\,(x,u)\in \cK^{c}\times\Act\,,\\[5pt]
\Lg^{u}\Lyap(x) &\,\le\,\kappa_2\bigl(1 + \rc(x,u)\bigr)
\qquad\forall\,(x,u)\in\cK\times\Act\,.
\end{split}
\end{equation}
\item[(iv)]
The functions $\upsigma$ and $\frac{\grad\Lyap}{1+\phi(\Lyap)}$
are bounded on $\Rd$.
\end{itemize}
\end{assumption}

Observe that when $\cK=\RR^{d}$ then the problem reduces to an
ergodic control problem with near-monotone running
cost, whereas if $\cK$ is bounded, we obtain an ergodic control problem for a
uniformly stable controlled diffusion.

\begin{theorem}
Let Assumption~\ref{A3.1} hold.
Then
\begin{itemize}
\item[(a)]
There exists $v^*\in\Ussm$ which is average cost optimal.
\item[(b)]
Every average cost optimal stationary Markov control is stable
and pathwise optimal.
\end{itemize}
\end{theorem}

%%%%%%%%%%%%%%%%%%%%%%%%%%%%%%%%%%
%\addtolength{\textheight}{-13.6cm}
%%%%%%%%%%%%%%%%%%%%%%%%%%%%%%%%%%

\begin{IEEEproof}
Part (a) follows by \cite[Theorem~3.1]{ABP14}, which
also shows that an average cost optimal stationary Markov control is
necessarily stable.

To prove part~(b), define $\psi_{N}\colon\RR_+\to\RR_+$,
for $N\in\NN$, by
\begin{equation*}
\psi_{N}(z) \,\df\, \int_{0}^{z} \frac{1}{N+\phi(y)}\,\D{y}\,,
\end{equation*}
and $\varphi_{N}\colon\Rd\to\RR_+$ by
\begin{equation*}
\varphi_{N}(x)\,\df\, \psi_{N}\bigl(\Lyap(x)\bigr)\,,\qquad x\in\Rd\,.
\end{equation*}

Let $U\in\Uadm$ be some admissible control such that
\begin{equation}\label{E-finite}
\liminf_{n\to\infty}\int_{\Rd\times\Act} \rc(x,u)\,\zeta^U_{t_n}(\D{x},\D{u})
\,<\,\infty\,,
\end{equation}
for some increasing divergent sequence $\{t_n\}$.
Since $\rc$ is inf-compact on $\cK\times\Act$,
it follows that $\zeta^U_{t_n}$ is a.s.\ tight
when restricted to $\Bor(\cK\times\Act)$.
Since $\grad\varphi_N$ and $\upsigma$ are bounded, by It\^o's formula
and \eqref{E-KA1} we obtain
\begin{align}\label{E-Ito2}
\frac{\varphi_N(X_t) - \varphi_N(X_0)}{t}
&\,=\, \frac{1}{t}\;\int_0^t \Lg^{U_s} \varphi_N(X_s)\,\D{s}\nonumber\\[5pt]
&\mspace{30mu}+
\frac{1}{t}\;\int_0^t \langle \grad\varphi_N(X_s),\upsigma(X_s)\, \D W_s\rangle\,.
\end{align}
Arguing as in \cite[Lemma~3.4.6]{book}, it follows that the second term on the
right hand side of \eqref{E-Ito}
converges a.s. to $0$ as $t\to\infty$.
We have
\begin{equation*}
\Lg^{u} \varphi_N(x) \,\le\, \begin{cases}
\frac{\kappa_0-\kappa_1\phi(\Lyap(x))}{N+\phi(\Lyap(x))}
+ h_{N}(x)\qquad\text{on~}\cK^c\,,\\[5pt]
\frac{\kappa_2(1 + \rc(x,u))}{N+\phi(\Lyap(x))}
+ h_{N}(x)\qquad\text{on~}\cK\,,
\end{cases}
\end{equation*}
where
\begin{equation*}
h_{N}(x)\,\df\,-\phi'\bigl(\Lyap(x)\bigr)\,
\frac{\babs{\upsigma\transp(x)\nabla\Lyap(x)}^2}{\bigl(N+\phi(\Lyap(x))\bigr)^2}\,,
\qquad x\in\Rd\,.
\end{equation*}
Hence, by \eqref{E-finite}--\eqref{E-Ito2}, for some constant $C$ we obtain
\begin{equation*}
\limsup_{n\to\infty}\;
\int_{\cK^c\times\Act} \frac{\phi\bigl(\Lyap(x)\bigr)}
{N+\phi\bigl(\Lyap(x)\bigr)}\,\zeta^U_{t_n}(\D{x},\D{u})
\,\le\, \frac{C}{N}\,,
\end{equation*}
from which it follows that $\zeta^U_{t_n}$ is a.s.\ tight
when restricted to $\Bor(\cK^c\times\Act)$.
Therefore, it is a.s.\ tight in  $\Pm(\Rd\times\Act)$.
Hence, almost surely
 every limit point of $\{\zeta^U_{t_n}\}$ is an ergodic occupation
measure \cite[Lemma~3.4.6]{book}.
Therefore
\begin{align*}
\liminf_{n\to\infty}\; \frac{1}{t_n}\int_{0}^{t_n} \rc(X_{s},U_{s})\,\D{s}
&\,\ge\, \inf_{\uppi\in\eom}\,\int_{\Rd\times\Act} \rc(x,u)\,\uppi(\D{x},\D{u})\\
&\,=\,\varrho_*\,.
\end{align*}
This completes the proof.
\end{IEEEproof}

\section{Conclusions}
It is well known that uniform stability, or equivalently,
uniform integrability of the hitting times, is equivalent to the existence
of a solution to the Lyapunov equation
\begin{equation}\label{E-Lyap4}
\Lg^u\Lyap(x) \,\le\, C- g(x)\qquad\forall\,(x,u)\in\Rd\times\Act\,,
\end{equation}
for some inf-compact function $g$ and a constant $C$.
Comparing \eqref{E-Lyap3} to \eqref{E-Lyap4}, it follows that
if $\Lyap$ grows no faster than $\E^g$ then the conclusions of
Corollary~\ref{C2.1} follow.
This closes significantly the gap between uniform stability
and sufficient conditions for pathwise optimality.

Even though the results of Section~\ref{S2} are specialized to controlled diffusions,
they are directly applicable to irreducible MDPs with
a countable state space, i.e., the model treated in \cite{Borkar-89}.
%%%%%%%%%%%%%%%%%%%%%%%%%%%%%%%%%%%%%%%%%%%%%%%%%%%%%%%%%%%%%%%%%%%%%%%%%%%%%%%%

\def\cprime{$'$} \def\cprime{$'$} \def\cprime{$'$}

\end{document}